\documentclass{article}

\usepackage{arxiv}

\usepackage[utf8]{inputenc} % allow utf-8 input
\usepackage[T1]{fontenc}    % use 8-bit T1 fonts
\usepackage{hyperref}       % hyperlinks
\usepackage{url}            % simple URL typesetting
\usepackage{booktabs}       % professional-quality tables
\usepackage{amsfonts}       % blackboard math symbols
\usepackage{nicefrac}       % compact symbols for 1/2, etc.
\usepackage{microtype}      % microtypography
\usepackage{lipsum}

\usepackage{graphicx} 
\usepackage{amsmath}
\usepackage{amssymb}
\usepackage{algorithm}
\usepackage{algorithmicx}
\usepackage{esint}% Include this line if your 

\if@secthm
  \newtheorem{thm}{Theorem}[section]
  \@addtoreset{thm}{section}
\else
  \newtheorem{thm}{Theorem}
\fi

\newtheorem{defn}[thm]{Definition}

  {\par
  \noindent
   {\bfseries\Elproofname}\enspace\ignorespaces}%
  {\par
  }
\def\Elproofname{PROOF.}

\title{An overview of SPSA: recent development and applications}

\author{
  Chen Wang
  Department of Mathematics\\
  University of Maryland\\
  College Park, MD 20046 \\
  \texttt{chen.wang@gmail.com} \\
}

\begin{document}
\maketitle

\begin{abstract}
There is an increasing need in solving high-dimensional optimization problems under non-deterministic environment. The simultaneous perturbation stochastic approximation (SPSA) algorithm has recently attracted considerable attention for solving high-dimensional optimization problems where the analytical formula cannot be attained. SPSA is designed to estimate the gradient by applying perturbation on a random subset of dimensions at each iteration. SPSA can be easily implemented and is highly efficient in that that it relies on measurements of the objective function, not on measurements of the gradient of the objective function. Since its invention, SPSA has been implemented in various fields such as reinforcement learning, production optimization etc. The paper briefly discuss the recent development of SPSA and its applications.
\end{abstract}


\section{Introduction}
The stochastic optimization problem setting can be described in the formula below 
\begin{equation}
\text{Find.  } \theta^\star = \arg min \mathbf{E}(f(\theta))
\end{equation}.

We are focusing on a setting in which the analytical form of the objective function of $f$ cannot be found but random
observations of the function can be obtained.
Meanwhile, this paper is concentrated on discussing the optimization algorithms where the gradient is involved, and we assume that noisy estimates of the objective function gradient are not directly available, so the function gradient needs to be estimated based on the noisy observations of the function.

The first two pioneers of stochastic approximation algorithms are the infinite small perturbation analysis developed by Robbins and Monroe \cite{robbins1951} and the finite difference method developed in \cite{kiefer1952}, also known as KW algorithm. Infinite-small perturbation analysis requires that the underline probability distribution of the noisy observations is known, whereas KW algorithms doesn't.

KW algorithm is usually more applicable since it requires less assumptions on the knowledge of the noisy observation $f(\theta)$. However, on the high dimensional cases, KW algorithm becomes very computationally expensive since it requires $2n$ rounds of simulation to get only one step of the gradient, with $n$ be the dimension of the parameter $\theta$. In order to accelerate the stochastic approximation process, many researchers have developed numerous approaches extending KW algorithm in hope of making the gradient estimation to be more efficient without sacrificing too much accuracy. A nominal example is the Simultaneous perturbation stochastic approximation(SPSA) developed by James Spall \cite{spall}.  Simultaneous perturbation stochastic approximation method is a random search technique that estimates the gradient using random perturbations that are independent, symmetric, zero-mean and satisfying an inverse moment bound. The most commonly used and studied class of perturbations within this category are those that are
independent, symmetric, ±1-valued, Bernoulli random variables. This algorithm (the standard SPSA, as
it is known), requires two function measurements at each update step.SPSA has become popular because of its
computational simplicity, as well as the convergence and rate guarantees that it provides.

Since the publish of SPSA algorithm, a lot of researchers have developed different extensions and alternatives. This paper will focusing on presenting the recent development on the approaches to tackle stochastic approximation on high-dimensional parameter.

The paper will be organized as follows: On section 2, we will discuss stochastic approximation methods with hessian estimation is involved. On section 3, we will discuss a SPSA alternative called random perturbation stochastic approximation. On section 4, we will discuss applying stocahstic approximation on a specific non-linear transform of random variable: cumulative prospect theory.

\section{Simultaneous perturbation stochastic approximation}
Let's denote the sample of $f(\theta)$ to be $\hat{f}(\theta)$, and assume $\theta$ is a $N$-dimensional vector. We also denote the N-dimensional noisy gradient of $\mathbf{E}(f(\theta))$ by $g(\theta)$.
In KW algorithm, noisy gradient estimator $\hat{g}(\theta)$ of $g(\theta)$ requires $2N$ noisy samples of $f$ with each dimension $\hat{g_i}(\theta)$ equals
\begin{equation}
    \hat{g}(\theta))_i = \frac{\hat{f}(\theta+c_ne_i)-\hat{f}(\theta-c_ne_i)}{2c_n}, 1 \leq i. \leq N
\end{equation}
Where $e_{i}$ is the unit vector with a 1 in the $i^{th}$ place, and $c_{n}$ is a small positive number that decreases with n.

It is apparently that with the dimension $N$ of $\theta$ getting large, KW algorithmic type estimator of the noisy gradient of $\mathbf{E}(f(\theta))$ is getting very inefficient. 

Simultaneous perturbation stochastic approximation method is achieved by introducing a perturbation vector $\Delta_n$ and apply it on each dimension of $\hat{g}(\theta)$. In another words, on SPSA setting, $\hat{g}(\theta)_i$ equals
\begin{equation}
\hat{g}(\theta))_i = \frac{\hat{f}(\theta+c_n\Delta_n)-\hat{f}(\theta-c_n\Delta_n)}{2c_n}.
\end{equation}
Note that SPSA essentially apply one numerator on all the dimensions of the quotient estimator of $\hat{g}(\theta)$. Henceforth, KW estimator perturbs only one direction at a time, while the SPSA estimator disturbs all directions at the same time (the numerator is identical in all p components). The number of loss function measurements needed in the SPSA method for each $\hat{g}(\theta)$ is always 2, independent of the dimension of theta $N$. Thus, SPSA uses $N$ times fewer function evaluations than KW estimator, which makes it a lot more efficient.

A good choice for $\Delta_{ni}$ is the Rademacher distribution, i.e. Bernoulli +-1 with probability 0.5. Other choices are possible too, but note that the uniform and normal distributions cannot be used because they do not satisfy the finite inverse moment conditions.

Since the SPSA gradient estimator can be proved to be an unbiased estimator of the gradient $g(\theta)$, the SPSA optimization algorithm can be proved to be converging to the optimal point of the function $\mathbf{E}(f(\theta))$ under certain differentiable conditions of $\mathbf{E}(f(\theta))$ . 

One example work of convergence proof of SPSA algoritm can be found in \cite{977290}. Mayak and Chen \cite{977290} denotes $\mathbf{J}(\theta)$ as $\mathbf{E}(f(\theta))$, it assume that $J(\theta)$ is twice differentiable and individual elements of the third derivative must be bounded. Additionally, $\nabla J$ must be Lipschitz continuous, bounded and the ODE $\dot{\theta}=g(\theta)$ must have a unique solution for each initial condition.
Under these conditions and a few others, SPSA algorithm will converge in probability to the set of global minima of $J(\theta)$

\section{Stochastic approximation with Hessian estimation included}
Recall the Hessian matrix of a function has the formula
\begin{align}
H f= \begin{bmatrix}
  \dfrac{\partial^2 f}{\partial x_1^2} & \dfrac{\partial^2 f}{\partial x_1\,\partial x_2} & \cdots & \dfrac{\partial^2 f}{\partial x_1\,\partial x_n} \\[2.2ex]
  \dfrac{\partial^2 f}{\partial x_2\,\partial x_1} & \dfrac{\partial^2 f}{\partial x_2^2} & \cdots & \dfrac{\partial^2 f}{\partial x_2\,\partial x_n} \\[2.2ex]
  \vdots & \vdots & \ddots & \vdots \\[2.2ex]
  \dfrac{\partial^2 f}{\partial x_n\,\partial x_1} & \dfrac{\partial^2 f}{\partial x_n\,\partial x_2} & \cdots & \dfrac{\partial^2 f}{\partial x_n^2}
\end{bmatrix}
\end{align}
in which each element represents the partial derivatives of two dimensions $i$ and $j$. The Hessian matrix has been first applied on optimization since Newton method, which is inspired by the second order Taylor expansion of a function. It is widely known that Newton type optimization involving Hessian matrix could converge to the optimal point much faster than simple gradient descent method when the underline function is relatively smooth. As a result, adaptive Newton-type schemes that estimate the Hessian using noisy objective function observations have also gathered considerable attention over the years.

The earliest such Hessian matrix estimation is developed by Fabian \cite{fabian1968}. Like the classical KW method, Fabian developed the the Hessian estimation using finite-difference method and required $O(2N
)$ samples of
the simulation at each update epoch. Spall \cite{880982} presented a simultaneous perturbation estimate of
the Hessian that was based on four noisy function simulation results. Among the four simulation results, two of them are also been used to estimate
the gradient. In the case when noisy gradient estimates are directly available, Spall \cite{880982} also presented a Newton-Ralphson scheme requiring three measurements.

\section{Stochastic optimization on Cumulative Prospect Theory functional}
Prospect theory was built by Kahneman and Tversky
\cite{Kahneman1979}. The modified version, cumulative prospect theory(CPT) is developed in early 1990s
\cite{Tversky1992}. There are mainly 3 assumptions of cumulative prospect theory(CPT): 1. It assumes that people tends to think of their possible outcomes relative to a certain value, rather than a final absolute value. 2. People tends to have different sensitivities towards gains and losses. In order words, the marginal impact of losses and gains are usually going with different direction. This phenomenon can also be described as ``risk-sensitive''. 3. People usually overweight the possibility of extreme events but underweight the possibility of high-chance event. In other words, people usually prospect the unlikely event to have a higher chance than what the probability of that event indicate. 

The 3 assumptions of CPT can be translated into the following mathematical expression:
1) The utility function has a reference point against which gains
and losses are evaluated, and this expression refers to the claim that people tend to think of possible outcomes usually relative to a certain reference point; 2) The utility function is concave on gains
and convex on losses, and this property just implies risk-avers behaviour of human being; 3) A probability
weighting function (cf. Def. \ref{def:probability-weighting-function})
that transforms the cumulative distribution function of a distribution such that
the probability of extreme is over weighted and the probability of common event is under weighted. We define the weighting function as the following \cite{LIN20181}:

\begin{defn}
\label{def:probability-weighting-function}A \textit{probability weighting
function}, $w$, is a monotonically non-decreasing continuous function
from $\left[0,1\right]$ to $[0,1]$ with $w\left(0\right)=0$ and
$w\left(1\right)=1$. 
\end{defn}

Let $X$ be a real random variable with a given probability distribution function, and $b$ to be denoted as a reference point together with two utilities functions $u^{+}$ and $u^{-}$. The $w^{+}$ and $w^{-}$ denote two different smooth probability weighting function. The CPT-functional applied on the random variable $X$ can be expressed as the following \cite{LIN20181}: 
\begin{flalign}
\mathbb{C}\left(X\right) & =\int_{0}^{\infty}w^{+}\left(\mathnormal{P}\left(u^{+}\left(\left(X-b\right)_{+}\right)>x\right)\right)dx\nonumber \\
 & -\int_{0}^{\infty}w^{-}\left(\mathnormal{P\left(u^{-}\left(\left(X-b\right)_{-}\right)>x\right)}\right)dx,\label{eq:CPT_main}
\end{flalign}
Notice that the notations $\left(\cdot\right)_{+}$
and $\left(\cdot\right)_{-}$ are shorthands for $\max\left(\cdot,0\right)$
and $\mbox{-\ensuremath{\min\left(\cdot,0\right)}}$, respectively.
Appropriate integrability assumptions are satisfied. 

In a lot of applications, for example financial markets and traffic light control system, one needs to find the best policy or strategy to optimize the system in accordance of CPT-functional so that human's satisfaction can be maximally achieved. That can be translated to a mathematical problem where the outcome of a system under difference policy can be modeled by a random variable $X(\theta)$, where $\theta$ is a parameters-based representation of a given policy. The optimization problem can be mathematically formulated by 
\begin{equation}
    \text{Find.  } \theta^\star = \arg min \mathbb{C}(X(\theta))
\end{equation}

Notice that unlike traditional expectation of a random variable, CPT functional is a non-linear transform of the underline random variable, therefore traditional Robins Monroe algorithm cannot be applied on deriving an analytical form of the CPT gradient. Since the only viable option is a finite difference variant of gradient estimator, \cite{8329994} derived a finite difference alike SPSA method on optimizing the CPT functional. Also notice that the estimating the gradient of the CPT functional using finite difference method requires an estimation of the CPT functional first, and \cite{la2016cumulative} firstly derived an estimator of the CPT functional based on quantile statistics. The proof of convergence of the CPT estimator can be established by the classical empirical statistics theory. Meanwhile, when implying certain assumptions on the CPT-functional like \cite{977290}, the paper \cite{8329994} proved that the established SPSA-alike gradient descent method can converge to the local minima of CPT-functional.

\bibliographystyle{unsrt}  
\bibliography{references} 

\end{document}